\documentclass{article}
\usepackage{heron2e}
\usepackage{times}

\title{Cohomological Properties of Differential
Calculi\\ on Hopf Algebras}  
\author{
F. Bonechi\adr{1,2}, R. Giachetti\adr{1}, R. Maciocco\adr{1}, 
E. Sorace\adr{1} and M.Tarlini\adr{1}}

\address[1]{Dipartimento di Fisica, Universit\`a di Firenze\\
             and I.N.F.N. Sezione di Firenze}
\address[2]{Laboratoire de Physique Th\'eorique,
Universit\'e de Paris VII, France}

\newfont{\ccc}{cmcsc10}

\newcommand{\id}{{\rm id}}
\newcommand{\gin}{\,{}_{\scriptstyle inv}\Gamma}
\newcommand{\pro}[1]{{\ccc (#1) Proposition.}} 
\newcommand{\cor}[1]{{\ccc (#1) Corollary.}} 
\newcommand{\lem}[1]{{\ccc (#1) Lemma.}}

\newcommand{\rmk}[1]{{\ccc (#1) Remark.}}

\newcommand{\dq}{{\cal D}}
\newcommand{\fq}{{\cal F}}
\newcommand{\ke}{{\rm Ker}\epsilon}
\newcommand{\uq}{{{\cal U}}}
\newcommand{\gd}{{}_\Gamma\hspace{-1pt} \delta}

\newcommand{\ad}{{\rm ad}}

\newcommand{\tf}{\tau_{\hspace{-1pt}\scriptscriptstyle {\cal F}}}
\newcommand{\tu}{\tau_{\scriptscriptstyle {\cal U}}}
\newcommand{\td}{\tau_{\hspace{-1pt}\scriptscriptstyle {\cal D}}}

\newcommand{\C}{{\cal C}}

\begin{document}
\maketitle
%\begin{abstract}
%  Abstracts are not required. We leave it to the author(s) to decide
%  whether to include an abstract or not.
%\end{abstract}

\section{Introduction}
In the approach to the differential calculus on quantum groups
proposed in \cite{swor1}, besides the obvious notion of differential $d$,
the theory is founded on the notion of bicovariant
bimodule; the algebraic nature of this construction extends the 
properties of differential forms to the noncommutative situation.
Following this idea a general treatment and a classification of
differential calculi have been constructed for quantum groups obtained
as deformation of semisimple Lie groups making use of the
quasitriangularity property \cite{bjur,ksch}.

In this report we give an intrinsic treatment of the
results we developed in {\cite{fbon} connecting the differential
calculi on Hopf algebras to the Drinfeld double \cite{vdri}.
In the first place we recover that bicovariant bimodules are in one
to one correspondence with the Drinfeld double representations
\cite{swor2}; we then introduce a Hochschild cohomology of the algebra
of functions and discuss the main result stating that each differential
calculus is associated to a 1-cocycle satisfying an additional
invariance condition with respect to a natural action \cite{fbon}.
Defining a Hochschild cohomology of the double, the above invariance
becomes a condition with respect to the enveloping algebra 
component of the double that must be added to the 1-cocycle
relation.  

The general classification of differential calculi is therefore 
reduced to a cohomological problem,
which can be performed with the more usual and efficient tools. 
Moreover a supply of differential calculi is obtained by 
observing that the coboundary operator maps invariant 0-cochains 
into invariant 1-coboundaries.

The construction we present is completely independent of the  
quasi-triangular property of Hopf algebras and can obviously be applied
to classical groups, both Lie and discrete or finite. An interesting
feature is that all the known differential calculi on quantum
and finite groups correspond to coboundaries, at difference with the
usual Lie group case, in which no invariant coboundary exists and the
classical differential calculus is determined by a nontrivial 1-cocycle.

\section{Differential calculus of the first order}
Let ${\cal A}$ be an associative algebra and $\Gamma$ an $\cal A$ -
bimodule; with $a\,.\,\gamma$ and $\gamma\,.\,a$ we
indicate left and right multiplication of $a\in{\cal A}$ with 
$\gamma\in\Gamma$.  

By {\it differential calculus of the first order} we mean the couple
($\Gamma,d$) where $d: {\cal A}\rightarrow \Gamma$ is a linear application
such that

$(i)\ $ $d(ab) = a\,.\,(db)+(da)\,.\,b \;\;\;a,b \,\in\,{\cal A}$; 

$(ii)$ Im $d$ generates $\Gamma$ ({\it i.e.} $\forall
\, \gamma \,\in\Gamma$ there exist $a_k, b_k \,\in\,{\cal A}$ such that
$\gamma = \sum_k a_k d b_k$).

Two differential calculi ($\Gamma_1,d_1$) and ($\Gamma_2, d_2$) are said
to be equivalent if there exists a bimodule isomorphism $i:\, \Gamma_1
\rightarrow \Gamma_2$ such that $d_2 = i \circ d_1$.

Let ${\cal A}_2 = {\rm Ker}\,m \subset {\cal A}\otimes {\cal A}$,
where $m:{\cal A}\otimes{\cal A}\rightarrow {\cal A}$ is the  multiplication 
in ${\cal A}$. Let  ${\cal A}_2$  be given the following structure
of $\cal A$ - bimodule:
\begin{eqnarray*}
a\,.\,(\sum_k x_k \otimes y_k)  &=& \sum_k a\, x_k 
\otimes y_k\ ,   \\
(\sum_k x_k \otimes y_k)\,.\,a &=& \sum_k x_k \otimes y_k\, a\ ,
\end{eqnarray*}
with $ \sum_k x_k y_k = 0$ and $a, b \in {\cal A}$.
It is easy to verify that, with 
$$ Da = 1 \otimes a - a \otimes 1\ . $$
(${\cal A}_2, D$) is a first order differential calculus.
We refer to it as to the {\it universal differential calculus}.
The universality property is justified by the following

\smallskip
\pro{1} {\it Each differential calculus {\rm (}$\Gamma,d${\rm )} is 
obtained as a quotient
from (${\cal A}_2, D$), {\rm i.e.} there exists a sub-bimodule
${\cal N} \subset {\cal A}_2$, such that $\Gamma = {\cal A}_2/{\cal N}$ and
$d = \pi \circ D$, where
$\pi:{\cal A}_2 \rightarrow \Gamma$ is the canonical homomorphism.}

\medskip
In the case of an Hopf algebra $\fq$ the universal calculus can be 
described in a particularly useful form. Let $r: \fq\otimes\fq 
\rightarrow \fq\otimes\fq$ be defined by
$$r( a\otimes b) = (a\otimes 1) \Delta b\ ,  $$     
with inverse
$$r^{-1} ( a\otimes b) = (a\otimes 1) (S \otimes 1) \Delta b  \ .   $$
Define on $\fq \otimes \ke$ the structure of $\fq$ - bimodule 
as follows
\begin{eqnarray*}
a\,.\,(\sum_k b_k \otimes c_k)  &=& \sum_k a\, b_k \otimes c_k\ , \\
(\sum_k b_k \otimes c_k)\,.\, a &=&  
\sum_k \left(b_k \otimes c_k \right)\Delta a\ ,
\end{eqnarray*}
with $\sum_k b_k \otimes c_k \, \in \, \fq \otimes \ke$. It's easy to prove
that $r$ is a bimodule isomorphism between $\fq_2$ and $\fq\otimes\ke$.
Moreover, if $D' = r\circ D : \fq \rightarrow \fq \otimes \ke\, ,$
{\it i.e.}
$$D'a = \Delta a - a\otimes 1  \;,$$
then ($\fq_2,D$) and ($\fq\otimes\ke,D'$) are isomorphic differential 
calculi. From now on we refer to ($\fq\otimes\ke,D'$) as to the 
universal differential calculus

\section{Bicovariant bimodules and Drinfeld quantum double}
The notion of {\it bicovariant bimodule} is an {\it ad hoc} notion
introduced in \cite{swor1} to generalize the definition of
translations on the space of forms on Lie groups to the case of
Hopf algebras. We show in this section the general result that 
bicovariant bimoduli are completely characterized by the
representations of the Drinfeld double of $\fq$.

\medskip

An $\fq$-bimodule is said to be a {\it left-covariant} bimodule if
it is defined a coaction 
$\delta_\Gamma : \Gamma \rightarrow \fq \otimes \Gamma$ such that
$$\delta_\Gamma(a\,\gamma)=
\Delta(a)\,\delta_\Gamma(\gamma)\,, ~~~~~~~\delta_\Gamma(\gamma\, a)=
\delta_\Gamma(\gamma)\,\Delta(a)\,,~~~~~~~
(\epsilon\otimes\id)\;\delta_\Gamma(\gamma)=\gamma\,,$$
for each $a\,\in\,\fq$ e $\gamma \,\in\,\Gamma$.
Two left-covariant $\fq$-bimoduli {\rm ($\Gamma$, $\delta_\Gamma$)} and 
{\rm ($\Gamma'$, $\delta_\Gamma'$)} are equivalent if there exists
a bimodule isomorphism $i: \Gamma\rightarrow \Gamma'$ such that
$$
(\id\otimes i) \delta_\Gamma = \delta_\Gamma' \, i  \ .
$$
A form $\omega\in\Gamma$ is left invariant if 
$\delta_\Gamma(\omega) = 1\otimes\omega$; we call $\gin$ the space
of left invariant forms.

Left covariant bimoduli are completely determined by giving a 
representation of $\fq$ on $\gin$, according to the following result 
(\cite{swor1}, but for the intrinsic formulation see \cite{cber}).

\smallskip
\pro{2} {\it 
If $\Gamma$ is a left covariant $\fq$-bimodule, then
$\tf: \fq \rightarrow {\rm End}(\gin)$ defined by  
$$
\tf (a) \gamma = \sum\limits_{(a)} S(a_{(1)}).\gamma.a_{(2)} \;\;\;
\gamma\,\in\gin\, , \,\, a \,\in\,\fq
$$
is right representation of $\fq$. Moreover
$\Gamma$ is isomorphic to
$\fq \otimes \gin$, which is a free left module and has a right
multiplication given by
$$
(1\otimes\gamma).a = 
\sum\limits_{(a)} a_{(1)} \otimes \tf (a_{(2)}) \gamma  \;\;\;
a\,\in\,\fq,\,\, \gamma\,\in\gin \;.   \eqno(3)
$$

Conversely, if $\tf$ is a right representation of $\fq$ on $V_\tau$, then
{\rm (3)} defines on the free left module $\Gamma = \fq \otimes V_\tau$ 
the structure of left covariant $\fq$-bimodule, with $\gin=V_\tau$.

Two left covariant $\fq$-bimoduli are equivalent if and only if the
corresponding representations of $\fq$ are equivalent.}

\medskip
Analogously an $\fq$-bimodule is {\it right-covariant} if it is
defined a coaction $\gd : \Gamma \rightarrow \Gamma \otimes \fq$ 
such that
$$\gd(a\,\gamma)=
\Delta(a)\,\gd(\gamma)\,, ~~~~~~~\gd(\gamma\, a)=
\gd(\gamma)\,\Delta(a)\,,~~~~~~~
(\id\otimes\epsilon)\;\gd(\gamma)=\gamma\,,$$ 
for each $a\,\in\,\fq$ e $\gamma \,\in\,\Gamma$. We omit the obvious 
definition of isomorphism of right covariant bimoduli.

\medskip
We finally arrive to the definition of bicovariant bimodule. 
An $\fq$-bimodule is {\it bicovariant} if it is left and right 
covariant and if left and right translations $\delta_\Gamma$ and 
$\gd$ verify the following compatibility relation
$$
(\id\otimes\gd)\delta_\Gamma \,=\, 
(\delta_\Gamma\otimes\id)\gd  \ .
$$
Two bicovariant bimoduli are isomorphic if there exists a bimodule isomorphism
which is an isomorphism of left and right covariant bimoduli.

\medskip
As in the case of left covariant bimoduli, there exists a one-to-one
correspondence between bicovariant bimoduli and representations of
an Hopf algebra $\dq$, the Drinfeld double of $\fq$. 

Let $\uq$ be the Hopf algebra dual to $\fq$. Then the Drinfeld double of
$\fq$ is the Hopf algebra $\dq$ defined by the following requirements:

({\it i}) $\dq =\fq\otimes\uq$ as vector space; 

({\it ii}) $\fq$ and $\uq^{\rm op}$ (the opposite Hopf algebra of $\uq$) are 
Hopf subalgebras of $\dq$; 

({\it iii}) for each $a\in\fq$ and $X\in\uq$ we have
$$
 aX = \sum\limits_{(a)(X)} X_{(2)}a_{(2)}<X_{(1)},S^{-1}(a_{(3)})>
<X_{(3)},a_{(1)}> \;.
 \eqno(4)
$$
In the following we will indicate with $\widetilde S$ the antipode in $\dq$,
{\it i.e.} ${\widetilde S}(a) = S(a)$ if $a\in\fq$ and ${\widetilde
S}(X)=S^{-1}(X)$  if $X\in\uq$.
We remark that $\dq$ is a quasitriangular Hopf algebra even if $\fq$
is not. The main result of this section is contained in the following 
proposition, which we state without proof.

\smallskip
\pro{5} {\it 
If {\rm(}$\Gamma,\gd,\delta_\Gamma${\rm)} is a bicovariant $\fq$-bimodule
then $\tu: \uq\rightarrow {\rm End}(\gin)$ defined by
$$
\tu (X) \gamma = (\id\otimes{\widetilde S}(X)) \gd(\gamma) \;\;\;
X\in\uq,\, \gamma\in\gin\,,  \eqno(6)
$$
is a right representation of $\uq$ on $\gin$ and 
$\td: \dq \rightarrow {\rm End}(\gin)$, defined by
$$
\td(Xa) = \tf (a) \tu (X) \;\;\; a\in\fq\hookrightarrow\dq,\,
             X\in\uq\hookrightarrow\dq\,,
$$
is a right representation of $\dq$.

Conversely, if $\td$ is a right representation on $\dq$ on $V_\tau$, 
$\tu =\td|_{_\uq}$ and $\Gamma=\fq\otimes V_\dq$ is the
left covariant bimodule associated to $\tf =\td|_{_\fq}$, then the
right translation $\gd$ defined in {\rm(6)} gives $\Gamma$ the
structure of a bicovariant bimodule.

Two representations of $\dq$ are equivalent if and only if the
corresponding bicovariant bimoduli are equivalent. }

\section{Bicovariant differential calculi and invariant 1-cocycles}
Let ($\Gamma,\delta_\Gamma,\gd$) be a bicovariant $\fq$-bimodule. 
A differential calculus ($\Gamma,d$) is said to be bicovariant if
$$\delta_\Gamma \circ d = \left( 1 \otimes d\right)\Delta\,,
\quad\quad  
\gd \circ d = \left( d \otimes 1 \right)\Delta  \;. $$
In the case of classical differential calculus these conditions express
the commutativity between differential and translations.
The requirement of bicovariance strongly restricts the set of 
admissible calculi.
An intrinsic property that in principle classifies them was already 
given in \cite{swor1}.

\smallskip
\pro{7} {\it Bicovariant differential calculi are in one-to-one correspondence
with right ideals ${\cal R}\subset\ke$ of $\fq$ which are {\rm ad}-invariant, 
{\rm i.e.} such that $\ad({\cal R}) \subset {\cal R}\otimes\fq$.} 

\medskip
In this section we propose an alternative approach in terms of
cocycles of an Hochschild cohomology, invariant with respect to a certain 
action of $\uq$. Let's recall the basic definitions of the Hochschild
cohomology.

\smallskip
Let $\cal A$ be a generic associative algebra and ${\cal C}^k({\cal
A},{\cal M})$ 
the set of $k$-linear applications ({\it k-cochains}) from
${\cal A}^k$ to an ${\cal A}$-bimodule
${\cal M}$, as usual $\C^o({\cal A},{\cal M})\equiv{\cal M}$ and we let 
$\delta: {\cal C}^k \rightarrow {\cal C}^{k+1}$ be
the coboundary operator defined by
$$
\delta\psi\,(\alpha_1,\,\alpha_2,\,\dots\,,\,\alpha_{k+1})=
\alpha_1.\psi(\alpha_2,\,\dots\,,\,\alpha_{k+1})\ + 
\eqno(8)
$$
$$
\sum\limits_{i=1}^{k}\,(-1)^i\,
\psi(\alpha_1,\dots\,,\,\alpha_i\,\alpha_{i+1},\dots\,,\,
\alpha_{k+1}) +
(-1)^{k+1}\,
\psi(\alpha_1,\dots\,,\,\alpha_k)\cdot\alpha_{k+1}\,. 
$$

${\cal C}({\cal A},{\cal M})$ is a complex and we define
{\it k-cocycles} and {\it k-coboundaries} as
$$Z^k({\cal A},{\cal M})= \{\psi \in {\cal C}^k \,|\, \delta \psi =
0\}\,,\quad\quad B^k({\cal A},{\cal M})=\{\psi \in {\cal C}^k \,|\, 
\exists \,\gamma\,\in\,{\cal C}^{k-1} \, | \,\psi 
\,=\, \delta \gamma\}$$
respectively. The $k$-th {\it group of cohomology} is
$$H^k({\cal A},{\cal M}) = Z^k/B^k .$$

Let $\fq$ be an Hopf algebra and ($\Gamma,d$) a bicovariant differential
calculus. Let's give to $\gin$ the structure of $\fq$ bimodule,
with right multiplication given by $\tf$ and left multiplication by
the counit $\epsilon$. Let's consider the Hochschild complex 
${\cal C}(\fq,\gin)$, in particular for $\gamma\in\gin$, we have
$(\delta\gamma)(a)=\epsilon(a) \gamma - \tf (a)\gamma  \;.$ 

Let $\bullet: {\cal C}^k(\fq,\gin)\otimes\uq\rightarrow 
{\cal C}^k(\fq,\gin)$ be the right action of $\uq$ on $k$-cochains 
defined by  
$$
(\psi\bullet X) (a_1,\ldots,a_k) = \sum_{(X)} \td (X_{(k+1)})
\psi(Ad_{\scriptscriptstyle X_{(k)}} a_1,\ldots,
Ad_{\scriptscriptstyle X_{(1)}} a_k ) \ .
$$

If $ \psi\bullet X = \epsilon (X) \psi$, for each $X\in\uq$, we say
that $\psi$ is {\it invariant}. 

\smallskip
\lem {9} {\it A $k$-cochain $\psi \in {\cal C}^k(\fq,\gin)$ is invariant
if and only if, for each $X\in\uq$,
$$
\td (S(X)) \psi(a_1,\ldots,a_k) = \sum_{(X)} 
\psi(Ad_{\scriptscriptstyle X_{(k)}} a_1,\ldots,
Ad_{\scriptscriptstyle X_{(1)}} a_k )  \;.
$$}

\medskip
We finally give the main result of this paper that completely
characterizes bicovariant differential calculi.

\smallskip
\pro{10} {\it
$(i)$ There is one-to-one correspondence between bicovariant differential
calculi {\rm ($\Gamma$,$d$)} and invariant $1$-cocycles $\psi$ such that
$$
da = \sum_{(a)} a_{(1)}.\psi(a_{(2)})  \; .
$$
$(ii)$ The image of invariant $0$-cochain under $\delta$ is just the
space of invariant $1$-cocycle, so that each invariant $0$-cochain
defines a {\rm coboundary differential calculus}.
}

\medskip
Using the whole representation of $\dq$ we can give $\gin$ the structure
of a $\dq$-bimodule, where $\alpha.\gamma = \epsilon(\alpha)\gamma$ and 
$\gamma.\alpha=\td (\alpha)\gamma$, for $\alpha\in\dq$ and 
$\gamma\in\Gamma$. The result in (10) can be equivalently 
reformulated in terms of Hochschild Cohomology $\C(\dq,\gin)$. 

\smallskip
\cor{11} {\it Bicovariant differential calculi are in one-to-one correspondence
with $\phi\in Z^1(\dq,\gin)$ such that $\phi(\uq)=0$.}

\medskip
\rmk{12} ({\it Coboundary calculi}) Two calculi that differ by a 
coboundary are not equivalent. Indeed consider an invariant 0-cochain
$\gamma$ and the calculus associated to $\delta\gamma$.
Using the property of invariance of $\gamma$,
$\td (X)\gamma = (\id\otimes S(X))\gd(\gamma) = \epsilon(X) \gamma$, 
we see that $\gd(\gamma)=\gamma\otimes 1  \;,$
so that an invariant zero cochain is a form which is right and left
invariant. The associated differential 
$$
da = \sum\limits_{(a)} a_{(1)}.(\delta\gamma)(a_{(2)}) =
a.\gamma -\sum\limits_{(a)} a_{(1)}.\td (a_{(2)})\gamma =
a.\gamma-\gamma.a \,,
$$
results in an internal derivation but it is obviously not equivalent 
to the trivial calculus $d=0$. In \cite{ksch} this kind of calculus is
called {\it internal} and it is shown that all bicovariant calculi on
quantum groups of the series $ABCD$ are coboundary. The same situation
is met with finite groups.
On the contrary, the ordinary differential calculus on Lie groups, being
the space of forms symmetrical, it can't be coboundary.

\section{Quantum Lie algebras}
Let ($\Gamma,d$) be a bicovariant differential calculus and suppose that
${\rm dim}\, \gin=n<\infty$ so that we can introduce a basis
$\{\omega_i\}\ i=1,n$. We then have
$$
[\td (a)]_{ji}=\langle f_{ij}, a\ , \rangle\quad\quad
[\td (S(X))]_{ij}=\langle X,R_{ij}\rangle  \;,
$$
with $\{f_{ij}\},\ X \in\,\uq$ and $\{R_{ij}\},\ a\in\,\fq$.
Obviously the usual properties of representative elements are satisfied:
$$\Delta (f_{ij}) = \sum_k f_{ik}\otimes f_{kj}\,,~~~~~\quad 
\epsilon\,(f_{ij})=\delta_{ij}\,,$$
$$\Delta (R_{ij}) = \sum_k R_{ik}\otimes R_{kj}\,,~~~~~\quad 
\epsilon\,(R_{ij})=\delta_{ij}\,.$$
Using the representation $\td$ the relation (4) reads 
$$
\sum_i (f_{li}*a)\, R_{si}  =  \sum_i R_{il} \, (a*f_{is})\;,
$$
where $f*a=(\id\otimes f)\Delta a$ and $a*f=(f\otimes \id)\Delta a$.
This is the bicovariance as in \cite{swor1}.

Again in the basis $\{\omega_i\}$ we find
$$
\omega_i.a = \sum_j (f_{ij}*a).\omega_j \quad\quad
\gd(\omega_i) = \sum_j \omega_j\otimes R_{ji}  \;.
$$

And the bicovariance of the bimodule implies that the numerical matrix
$$
\Lambda^{ij}_{k\ell}=\langle f_{j\ell}\,,\,R_{ki}\rangle 
$$
solves quantum Yang Baxter equation
$\Lambda_{12}\Lambda_{13}\Lambda_{23}=\Lambda_{23}\Lambda_{13}\Lambda_{12}$.
We used the convention $(A\otimes B)_{j\ell}^{ik}=A_{ij} B_{k\ell}$. 
This is an obvious consequence of the quasitriangularity of $\dq$.
In fact if $\{e_A\}$ is a linear basis for $\fq$, $\{e^B\}$ for $\uq$, then 
${\cal R} = \sum_A e_A\otimes e^A \in \dq\otimes\dq$ and 
$\sigma\circ{\cal R}^{-1} = \sum_A S(e^A)\otimes e_A$, with 
$\sigma(a\otimes b)=b\otimes a$, solve the universal quantum Yang-Baxter 
equation. Using this basis it is possible to write
$$
R_{ij} = \sum_A e_A [\td (S(e^A))]_{ij} \quad\quad
f_{ij} = \sum_A e^A [\td (e_A)]_{ji}  \;,
$$
and then, introducing the left representation $\rho_\dq = \td^t$,
$$
 \Lambda_{k\ell}^{ij} = \rho_\dq[\sigma\circ{\cal R}^{-1}]^{ij}_{k\ell}  \;.
$$

\medskip
Let $\psi$ be the invariant 1-cocycle associated to 
($\Gamma,d$) and let $\chi_i\in\uq$ be defined by 
$$
\langle \chi_i,a \rangle = [\psi(a) ]_i \ .
$$
The cocycle properties of $\psi$ yield the relations
$$
(i)\quad \Delta(\chi_i) = 1\otimes\chi_i + \sum_j \chi_j \otimes f_{ji} \,,
\quad\quad\quad
(ii)\quad \epsilon(\chi_i) = 0  \ .
$$
The invariance of $\psi$ gives then 
$$
(iii)\quad ad_Y(\chi_i) = \sum\limits_{(Y)} S(Y_{(1)}) \chi_i Y_{(2)} = 
\langle Y,R_{ik}\rangle \,\chi_k\;,   \quad Y\in\uq \;.
$$
that constitutes an invariant condition of the ``{\it fields}'' $\chi_i$    
with respect to the adjoint action of $\uq$ on itself.
Properties ({\it i-iii}) identify what is usually called a 
{\it quantum Lie algebra} (see \cite{dber}). 

The invariants fields
allow to the following expression of the differential:
$$ da = \sum_i (\chi_i * a).\omega_i    \ .$$

\section{Differential calculus for the finite groups}
The results of the previous sections are well illustrated on finite
groups. In fact it is possible to give a complete classification
of bicovariant bimoduli by means of representations of the double; it is
also possible to study the invariant cocycles of the associated 
Hochschild cohomology, and then obtain a complete classification
of bicovariant differential calculi of finite groups.

Let $G$ be a finite group of order $n$ and $\fq(G)$ the $n$-dimensional
vector space of complex functions on $G$. The elements 
$\{\phi_g\}_{g\in G}$, where $\phi_g(h)=\delta_{g,h}$, are a basis for
$\fq(G)$. 
As usual the Hopf algebra structure on $\fq(G)$ is given by
$$
\phi_g \phi_h = \delta_{g,h} \phi_g\, , \quad 
\Delta \phi(g,h) = \phi(gh)\, , \quad S(\phi)(g) = \phi(g^{-1})\,, 
\quad \epsilon(\phi) = \phi(e)\,,
$$
where $g,h \in G$ and $e$ is the unity of the group. 

We can canonically associate to $\fq(G)$ another Hopf algebra, the group
algebra {\bf C}$G$. The multiplication is induced by the group
composition and the coalgebra structure is determined by
$$
\Delta g = g\otimes g\ ,\quad S(g) = g^{-1}\ , \quad \epsilon(g)=1,   
\quad\quad\quad g\in G \,.
$$ 
 {\bf C}$G$ and $\fq(G)$ are in duality with the coupling 
$$
\langle g,\phi \rangle = \phi(g)  \quad  \phi\in\fq(G)\ , \quad g\in G  \;.
$$
In particular $\{\phi_g\}_{g\in G}$ and $\{g\}_{g\in G}$ are dual basis,
{\it i.e.} $\langle g,\phi_h\rangle = \delta_{g,h}\;.$
\smallskip
The double $\dq(G)$ is isomorphic to $\fq(G)\otimes {\bf C}G$ as a vector
space and contains $\fq(G)$ and {\bf C}$G$ as Hopf subalgebras. Its 
definition is completed by the relations
$$
g \phi = \phi^g g  \quad\quad   \phi\in\fq(G), \, g\in {\bf C}G  \;,
$$
where $\phi^g(h) = \phi(g^{-1} h g)$, with $g,h \in G$.

\smallskip
Let $\C = \{h\in G \,|\, h=ghg^{-1}\}$ be a conjugacy class and let 
$n_\C$ be the number of elements of $\C$. Let $a\in\C$ and $Z_a=\{g\in
G|ag=ga\}$ be its centralizer, the following holds

\smallskip
\pro{13} {\it Let $V_\C$ be a $n_\C$-dimensional vector space and 
$\{v_h\}_{h\in\C}$ a basis. 
For each irreducible representation of $G$ of dimension $n_\mu$,  
$\rho^\mu: G\rightarrow {\rm End}(W^\mu)$, an irreducible
representation of the double  
$\rho^\mu_\C: \dq(G)\rightarrow {\rm End}(V_\C\otimes W^\mu)$ 
is defined by
\begin{eqnarray*}
\rho^\mu_\C(\phi) \,v_h\otimes w_\alpha &=& \phi(h)\, v_h \otimes 
w_\alpha \\
\rho^\mu_\C(g) \, v_h\otimes w_\alpha &=&  v_{ghg^{-1}} \otimes 
\rho^\mu(g) \,w_\alpha \;,\quad\quad \alpha=1,\ldots,n_\mu\ .
\end{eqnarray*}
Two representations $\rho^\mu_{\C_1}$ and $\rho^\nu_{\C_2}$ 
are equivalent if and only if $\C_1=\C_2$ and 
the restrictions $\rho^\mu|_{Z_a}$ and $\rho^\nu|_{Z_a}$ are 
equivalent for an arbitrary $a\in\C$.}

\medskip
Using propositions (2) and (5) it is now easy to classify
bicovariant bimoduli. Let $\Gamma^\mu_\C$ the bicovariant bimodule
associated to $\rho^\mu_\C$; it is generated as a free left module
by left invariant forms $\{\omega_{[h\alpha]}\}$, $h\in\C$ and 
$\alpha=1,\ldots,n_\mu$. Right multiplication and right translation are
defined by
\begin{eqnarray*}
 \omega_{[h\alpha]}.\phi &=& \sum\limits_{\ell\in\C, \beta}
\left(F^{[\ell\beta]}_{[h\alpha]}*\phi\right).\omega_{[\ell\beta]} \\
 \gd(\omega_{[h\alpha]}) &=& \sum\limits_{\ell\in\C, \beta}
 \omega_{[\ell\beta]} \otimes R^{[\ell\beta]}_{[h\alpha]}\;,
\end{eqnarray*}
where
$$
F^{[\ell\beta]}_{[h\alpha]} = h
\,\delta_{h,\ell}\,\delta_{\alpha,\beta}\ ,\quad\quad\quad
R^{[\ell\beta]}_{[h\alpha]} = \sum\limits_{g\in G} 
\phi_g \, \delta_{\ell,ghg^{-1}} [\rho^\mu(g)]^\beta_\alpha  \ .
$$
Finally we conclude that the $(n_\mu n_\C)\times(n_\mu n_\C)$ matrix
$$
\Lambda^{[n\eta][\ell\beta]}_{[m\gamma][k\alpha]}\ =\  
\langle F^{[\ell\beta]}_{[k\alpha]}, R^{[m\gamma]}_{[n\eta]} \rangle  
\ =\  \delta_{\ell,k} \, \delta_{\alpha,\beta} \, \delta_{m,knk^{-1}}\, 
[\rho^\mu(k)]^\gamma_\eta
$$
satisfies the quantum Yang-Baxter equation.

The study of the corresponding Hochschild cohomologies permits the complete
classification of bicovariant differentials calculi on finite groups. Using 
proposition (10) in fact we arrive at the following result.

\smallskip
\pro{14} {\it
Let $\C$ a conjugacy class and $\Gamma_\C$ the bicovariant bimodule
associated to $\C$ and to the trivial representation of $G$.
Then {\rm($\Gamma_\C,d_\C$)}, with $d_\C:\fq(G)\rightarrow\Gamma_\C$ 
defined by 
$$
d_\C \phi = \sum\limits_{g\in\C} (\chi_g*\phi).\omega_g 
\quad\quad   \phi\in\fq(G) \;,
$$
where $\chi_g=e-g$, is a coboundary bicovariant differential calculus,
associated to the {\rm0}-cochain $\sum_g \omega_g$.
All non equivalent bicovariant differential calculi on $\fq(G)$ are direct
sums of these calculi.}

\end{document}